\documentclass[12pt]{article}
\usepackage{amssymb,amsfonts,amsmath, psfrag,cite}
\usepackage{graphicx, lscape}

\newtheorem{theo}{Theorem}

\makeatletter \@addtoreset{equation}{section}
\@addtoreset{theo}{section} \makeatother

\def\qed{\hfill \rule{4pt}{7pt}}
\def\pf{\noindent {\it Proof.} }

\def\C{{\mathbb C}}
\def\Z{{\mathbb Z}}
\def\N{{\mathbb N}}

\begin{document}

\begin{center}
{\Large\bf  Nonterminating  Basic Hypergeometric Series

and the $q$-Zeilberger Algorithm}

\vskip 3mm

 William Y.C. Chen$^1$,
Qing-Hu Hou$^2$ and Yan-Ping Mu$^3$\\[5pt]
Center for Combinatorics, LPMC \\
Nankai University, Tianjin 300071, P. R. China

\vskip 3mm

 E-mail: $^1$chen@nankai.edu.cn,
$^2$hou@nankai.edu.cn, $^3$myphb@eyou.com

\end{center}

\begin{abstract}
We present a systematic method for proving  nonterminating basic
hypergeometric identities. Assume that $k$ is the summation index.
By setting a parameter $x$ to $xq^n$, we may find a recurrence
relation of the summation by using the
 $q$-Zeilberger algorithm. This method applies to
  almost all nonterminating basic
hypergeometric summation formulas  in the book of Gasper and
Rahman. Furthermore, by comparing the recursions and the limit
values, we may verify many classical transformation formulas,
including  the Sears-Carlitz transformation, transformations of
the very-well-poised $_8\phi_7$ series, the Rogers-Fine identity,
and the limiting case of Watson's formula that implies the
Rogers-Ramanujan identities.
\end{abstract}



\section{Introduction and Notations}

This paper is devoted to develop a systematic method for proving
nonterminating basic hypergeometric series summation and
transformation formulas. The idea of finding recurrence relations
and proving basic hypergeometric series identities by iteration
has been used very often, see \cite{Andr96,Andr96b,Andr98,Fine}.
However, there does not seem to exist a systematic method within
the scope of computer algebra for proving nonterminating
hypergeometric summation and transformation formulas. One obstacle
lies in the infinity of the summation ranges. In this paper, we
find that the $q$-Zeilberger algorithm can be used as a mechanism
for proving many basic hypergeometric summation and transformation
formulas.

 Wilf and
Zeilberger developed an algorithmic proof theory for identities on
hypergeometric series and basic hypergeometric series
\cite{Zeil91,WZ92a,PWZ96}. For the purpose of this paper, we are
concerned with the $q$-Zeilberger algorithm. Koornwinder
\cite{Koor93}, Paule and Riese \cite{PauRie95}, and B\"{o}ing and
Koepf \cite{BoiKoe99} further studied algorithmic proofs of basic
hypergeometric identities. Most of the theory and implementations
are restricted to the case of terminating identities. For
nonterminating identities, Gessel \cite{Ges95} and Koornwinder
\cite{Koor98} provided computer proofs of Gauss' summation formula
and Saalsch\"{u}tz' summation formula by means of a combination of
Zeilberger's algorithm and asymptotic estimates.

Our approach can be described as follows. Let
\[ f(a,\ldots,c)=\sum\limits_{k=0}^\infty t_k(a,\ldots,c)\]
 be a
hypergeometric series with parameters $a,\ldots,c$. We aim to find
a recurrence relation of the form
\begin{multline} \label{rec}
p_0(a,\ldots,c) f(a,\ldots,c) + p_1(a,\ldots,c)
f(aq,\ldots,cq) + \cdots \\
+ p_d(a,\ldots,c) f(aq^{d},\ldots,cq^{d})  = 0,
\end{multline}
where $d$ is a positive integer and $p_0,\ldots,p_d$ are
polynomials. To this end, we try to find polynomials
$p_0,\ldots,p_d$ and a sequence $(g_0,g_1,\ldots)$ such that
\begin{multline}\label{t-g}
p_0(a,\ldots,c) t_k(a,\ldots,c) + p_1(a,\ldots,c)
t_k(aq,\ldots,cq) + \cdots \\
+ p_d(a,\ldots,c) t_k(aq^{d},\ldots,cq^{d})  = g_{k+1}-g_k.
\end{multline}
Assume that $g_0=\lim\limits_{k \to \infty} g_k = 0$. Then
\eqref{rec} follows immediately by summing over $k$ in
\eqref{t-g}.

The main idea of this paper is to use the $q$-Zeilberger algorithm
\cite{WZ92a,Koor93,PauRie95,BoiKoe99} to solve $p_i$ and $g_k$.
The bridge to the $q$-Zeilberger algorithm is the introduction of
a new variable $n$ by setting the parameters $a,\ldots,c$ to
$aq^n,\ldots,cq^n$. Then the summand $t_k(aq^n,\ldots,cq^n)$
becomes a bivariate $q$-hypergeometric term. Applying the
$q$-Zeilberger algorithm, we can always obtain an equation of form
\eqref{t-g}.

When the recurrence relation \eqref{rec} is of first order (i.e.,
$d=1$) or involves only two terms, $f(a,\ldots,c)$ equals its
limit value $\lim\limits_{N \to \infty} f(aq^N,\ldots,cq^N)$
multiplied by an infinite product. Using this approach, we can
derive almost all nonterminating summation formulas listed in the
appendix of the book of Gasper and Rahman \cite{GaRa90}, including
bilateral series formulas.

When the recurrence relation involves at least three terms, we
show that $f(a,\ldots,c)$ is determined uniquely by the recurrence
relation and its limit value under suitable convergence conditions
(Theorem~\ref{th-rec}). Therefore, to prove an identity, it
suffices to verify that both sides satisfy the same recurrence
relation and that the identity holds for the limiting case. Using
this method, we can prove many classical transformation formulas.

Let us  introduce some basic notation. The set of integers and
nonnegative integers are denoted by $\Z$ and $\N$, respectively.
Throughout the paper, $q$ is a fixed nonzero complex number with
$|q|<1$.

 The
$q$-shifted factorial is defined for any complex parameter $a$ by
\[ (a;q)_\infty = \prod_{k=0}^\infty (1-aq^k), \quad \mbox{and}
\quad (a;q)_n = {(a;q)_\infty \over (aq^n;q)_\infty}, \quad
\forall\, n \in \Z.
\]
For brevity, we write
\[
(a_1,\ldots,a_m;q)_n = (a_1;q)_n \cdots (a_m;q)_n,
\]
where $n$ is an integer or infinity. Furthermore, the basic
hypergeometric series are defined by
\[
{_r\phi_s} \left[ \begin{array}{c} a_1,\ldots,a_r\\
b_1,\ldots,b_s
\end{array}; q, z \right] =
\sum_{k=0}^\infty {(a_1,\ldots,a_r;q)_k \over
(b_1,\ldots,b_s;q)_k} {z^k \over (q;q)_k} \left( (-1)^k q^{k
\choose 2} \right)^{s-r+1},
\]
and the bilateral basic hypergeometric series are defined by
\[
{_r\psi_s} \left[ \begin{array}{c} a_1,\ldots,a_r\\
b_1,\ldots,b_s
\end{array}; q, z \right] =
\sum_{k=-\infty}^\infty {(a_1,\ldots,a_r;q)_k \over
(b_1,\ldots,b_s;q)_k} z^k  \left( (-1)^k q^{k \choose 2}
\right)^{s-r}.
\]

\section{Summation Formulas} \label{S2}

In this section, we present a method of proving nonterminating
basic hypergeometric identities by using the $q$-Zeilberger
algorithm. Given a term of the form
\[
{(a'_1,\ldots,a'_u;q)_\infty \over (b'_1,\ldots,b'_v;q)_\infty}
{(a_1,\ldots,a_r;q)_k \over (b_1,\ldots,b_s;q)_k} q^{d{k \choose
2}} z^k,
\]
by setting some parameters $a,\ldots,c$ to $aq^n,\ldots,cq^n$, we
get a bivariate $q$-hypergeometric term $t_k(aq^n,\ldots,cq^n)$ on
$n$ and $k$. By the $q$-Zeilberger algorithm, we obtain a
bivariate $q$-hypergeometric term $g_{n,k}$ and polynomials
$p_i(q^n,a,\ldots,c)$ which are independent of $k$
 such that
\begin{multline} \label{eq-ag}
p_0(q^n,a,\ldots,c) t_k(aq^n,\ldots,cq^n) + p_1(q^n,a,\ldots,c)
t_k(aq^{n+1},\ldots,cq^{n+1}) + \cdots \\
+ p_d(q^n,a,\ldots,c) t_k(aq^{n+d},\ldots,cq^{n+d})  =
g_{n,k+1}-g_{n,k}.
\end{multline}
Suppose that $g_{0,0} = \lim\limits_{k \to \infty} g_{0,k}=0$. By
setting $n=0$ in \eqref{eq-ag} and summing over $k$, we derive a
recurrence relation of form \eqref{rec} for
$f(a,\ldots,c)=\sum\limits_{k=0}^\infty t_k(a,\ldots,c)$.

When the recursion \eqref{rec} involves only two terms, say
$f(a,\ldots,c)$ and $f(aq^d,\ldots,cq^d)$, we have
\[
f(a,\ldots,c) = \lim_{N \to \infty} f(aq^{dN},\ldots,cq^{dN})
\cdot \lim_{N \to \infty} \prod_{i=0}^{N-1} \left(
-{p_d(aq^{di},\ldots,cq^{di}) \over p_0(aq^{di},\ldots,cq^{di})}
\right).
\]
Therefore, the evaluation of $f(a,\ldots,c)$ becomes the
evaluation of its limit value $\lim\limits_{N \to \infty}
f(aq^{dN},\ldots,cq^{dN})$, which is much simpler and is usually
 an infinite product.

\subsection{Unilateral Summations}

We now present an example to show how to obtain an infinite
product expression from an infinite summation.

{\noindent \bf Example.} The $q$-binomial theorem:
\[
f(a,z) = {_1\phi_0} \left[ \begin{array}{c} a\\
  \text{---}
  \end{array}; q, z \right] = {(az;q)_\infty \over (z;q)_\infty},
  \quad |z|<1.
\]

The $q$-binomial theorem was derived by Cauchy \cite{Cauchy},
Jacobi \cite{Jacobi} and Heine \cite{Heine}. Heine's proof
consists of using series manipulations to derive the recurrence
relation
\begin{equation}
\label{q-bino} (1-z)f(a,z) = (1-az) f(a,qz).
\end{equation}
Gasper \cite{Gasper} provided another proof using a recurrence
relation respect to the parameter $a$:
\[
f(a,z) = (1-az) f(aq,z).
\]

Our computer generated proof is similar to Heine's proof. The
recurrence relation generated by the $q$-Zeilberger algorithm
turns out to be  \eqref{q-bino}. Let $u_k(z)$ be the summand
\[
u_k(z) = {(a;q)_k \over (q;q)_k} z^k
\]
and $u_{n,k}=u_k(zq^n)$. By the $q$-Zeilberger algorithm, we
obtain
\begin{equation}
\label{qZ-bino}
(-azq^n+1)u_{n+1,k}+(zq^n-1)u_{n,k} = g_{n,k+1} -
g_{n,k},
\end{equation}
where $g_{n,k} = (1-q^k) u_{n,k}$. Denote the left hand side of
\eqref{qZ-bino} by $t_{n,k}$. Then
\[
\sum_{k=0}^\infty t_{n,k} = -g_{n,0} + \lim_{k \to \infty} g_{n,k}
= 0, \quad \forall\, n \ge 0,
\]
implying that
\[
f(a,zq^{n}) = {1-azq^n \over 1-zq^n} f(a,zq^{n+1}), \quad
\forall\, n \ge 0.
\]
Hence
\begin{eqnarray*}
f(a,z) &=& {1-az \over 1-z} f(a,zq) = {1-az \over 1-z} {1-azq
\over 1-zq}
f(a,zq^2) = \cdots \\
&=& \lim_{N \to \infty} {(az;q)_N \over (z;q)_N} f(a,zq^N) =
{(az;q)_\infty \over (z;q)_\infty},
\end{eqnarray*}
where $\lim\limits_{N \to \infty} f(a,zq^N) = 1$ holds by
Tannery's Theorem (see \cite[p.~292]{Tan04} or \cite{Boas65}):

\noindent {\bf Tannery's Theorem}. Suppose $s(n)=\sum_{k \ge 0}
f_k(n)$ is a convergent series for each $n$. If there exists a
convergent series $\sum_{k \ge 0} M_k$ such that $|f_k(n)| \le
M_k$, then
\[
\lim_{n \to \infty} s(n) = \sum_{k=0}^\infty \lim_{n \to \infty}
f_k(n).
\]

The following summation formulas (most of which come from the
appendix of \cite{GaRa90}) can be obtained by the  above method.
We only list the recursions and the limit values $\lim\limits_{N
\to \infty} f(aq^{dN},\ldots,cq^{dN})$.

\[
\begin{array}{|c|c|c|c|} \hline
\mbox{Name} & \mbox{Summation} & \mbox{Recurrence Relation} &
\parbox{40pt}{\centering Limit value}
  \rule{0pt}{15pt} \\[5pt] \hline
\mbox{$q$-exponential} & \sum\limits_{k=0}^\infty {z^k \over
(q;q)_k} &
  f(z)={1 \over 1-z} f(zq) & 1 \rule{0pt}{25pt} \\[10pt] \hline
\mbox{$q$-exponential} & \sum\limits_{k=0}^\infty {q^{k
 \choose 2} z^k \over (q;q)_k} &
  f(z)=(1+z)f(zq) & 1 \rule{0pt}{25pt} \\[10pt] \hline
\parbox{70pt}{\centering Lebesgue \\ \cite[p.~21]{Andrews76}} & \sum\limits_{k=0}^\infty {(x;q)_k q^{k+1 \choose
2}
  \over (q;q)_k} & f(x) = (1-xq)f(xq^2) & (-q;q)_\infty \rule{0pt}{25pt} \\[10pt] \hline
\parbox{70pt}{Generalization of Lebesgue   \cite{Andrews73,ChenLiu2}} & \sum\limits_{k=0}^\infty
{(a,b;q)_k
  q^{k+1 \choose 2} \over (q;q)_k (abq;q^2)_k} &
   \begin{array}{l} f(a,b) = \\[5pt]
    \quad {(1-aq)(1-bq) \over (1-abq)(1-abq^3)} f(aq^2,bq^2) \end{array}  & (-q;q)_\infty
  \rule{0pt}{30pt} \\[20pt] \hline
\mbox{$_1\phi_1$} & {_1\phi_1} \left[ \begin{array}{c} a\\ c
  \end{array}; q, \displaystyle {c \over a} \right] & f(c) = {1-c/a \over 1-c}
  f(cq) & 1 \rule{0pt}{25pt} \\[10pt] \hline
\mbox{$q$-Gauss} & {_2\phi_1} \left[ \begin{array}{c} a,b\\ c
  \end{array}; q, \displaystyle {c \over ab} \right] &
  f(c) = {(1-c/a)(1-c/b) \over (1-c)(1-c/ab)} f(cq) & 1
  \rule{0pt}{25pt} \\[10pt] \hline
\parbox{70pt}{ $q$-Kummer (Bailey-Daum) sum} &
  {_2\phi_1} \left[ \begin{array}{c} a,b\\ aq/b
  \end{array}; q, \displaystyle -{q \over b} \right] &
  \begin{array}{l} f(a) = \\[5pt]
  \quad {(1-aq^2/b^2)(1-aq) \over (1-aq^2/b)(1-aq/b)} f(aq^2) \end{array} &
  {(-q;q)_\infty \over (-q/b;q)_\infty}
  \rule{0pt}{30pt} \\[20pt] \hline
\parbox{70pt}{ A $q$-analogue of Bailey's $_2F_1(-1)$} &
  {_2\phi_2}\left[ \begin{array}{c} a,q/a\\ -q,b
  \end{array}; q, -b \right] &
  f(b) = {(1-ab)(1-bq/a) \over (1-bq)(1-b)} f(bq^2) & 1
  \rule{0pt}{25pt} \\[10pt] \hline
\parbox{70pt}{ A $q$-analogue of Gauss' $_2F_1(-1)$} &
  \begin{array}{r} {_2\phi_2} \scriptsize \left[ \begin{array}{c} a^2,b^2\\
  abq^{1 \over 2}, -abq^{1 \over 2}
  \end{array}; \right. \quad \\[5pt] \left. q, -q \right] \end{array} &
  \begin{array}{l} f(a,b) = \\[5pt]
  \quad {(1-a^2 q)(1-b^2 q) \over (1-a^2 b^2 q)(1-a^2 b^2 q^3)}
  f(aq,bq) \end{array} & (-q;q)_\infty
  \rule{0pt}{35pt} \\[20pt] \hline
\mbox{$q$-Dixon sum} & \begin{array}{r} \scriptsize {_4\phi_3} \left[ \begin{array}{c} a^2,-qa,b,c \\
  -a,a^2q/b,a^2q/c
  \end{array} ; \right. \quad \\[5pt] \left. q, {qa \over bc} \right]
  \end{array} &
  \begin{array}{l} f(a)  = \\[5pt]
  \quad {(1-a^2q/bc)(1-a^2q^2/bc)(1+aq) \over (1-a^2q/b)(1-a^2q^2/b)(1-aq/bc)}
  \\[5pt]
   \quad \ {(1-aq/b)(1-aq/c)(1-a^2q) \over
  (1-a^2q/c)(1-a^2q^2/c)}f(aq)
  \end{array}  & 1
  \rule{0pt}{40pt} \\[30pt] \hline
\end{array}
\]

Here are more examples.

\noindent {\bf A $q$-analogue of Watson's $_3F_2$ sum:}
\[
 f(a,c) = {_8\phi_7} \left[ \begin{array}{c} \mu, q
      \mu^{\frac{1}{2}}, -q \mu^{\frac{1}{2}},
      a^2,b^2q, c,-c,-abq/c \\
      \mu^{\frac{1}{2}},-\mu^{\frac{1}{2}},-bcq/a,-ac/b,-abq, abq, c^2
      \end{array}
      ; q, {c \over ab} \right],
\]
where $\mu = -abc$.

By the $q$-Kummer sum, we have
\[
\lim_{N \to \infty} f(aq^N,cq^N)= {_2\phi_1} \left[
\begin{array}{c} b^2q, -abq/c \\ -bcq/a \end{array}; q, {c \over
ab} \right] =
  {(-q;q)_\infty (b^2q^2,c^2q/a^2;q^2)_\infty \over
  (c/ab,-bcq/a;q)_\infty}.
\]
By computation, one derives that
\[ f(a,c) = {(1+abcq)(1+c/b)(1+abcq^2)(1-a^2q)(1-c/b) \over
      (1-c^2q)(1-abq)(1+abq)(1+acq/b)(1+ac/b)} f(aq,cq).
\]
Thus, we have
\[
f(a,c) = {(-abcq,-c/b,c/b,-q;q)_\infty
(a^2q,b^2q^2,c^2q/a^2;q^2)_\infty \over
(abq,-abq,-ac/b,c/ab,-bcq/a;q)_\infty (c^2q;q^2)_\infty}.
\]

\noindent {\bf A $q$-analogue of Whipple's $_3F_2$ sum:}
\[
f(c) =  {_8\phi_7} \left[ \begin{array}{c} -c,q(-c)^{\frac{1}{2}},
      -q(-c)^{\frac{1}{2}},a,q/a, c,-d,-q/d \\
      (-c)^{\frac{1}{2}},-(-c)^{\frac{1}{2}},-cq/a,-ac
      -q,cq/d,cd
      \end{array} ; q, c \right].
\]

By computation we have the recurrence relation
\begin{multline*}
f(c) = {(1+cq^2)(1+c)(1-cq^2/ad)(1-acq/d) \over
(1+acq)(1-cdq)(1+cq^2/a)(1-cq^2/d)} \\
\times {(1-acd)(1-cdq/a)(1+cq)^2
   \over (1+ac)(1-cd)(1+cq/a)(1-cq/d)}
   f(cq^2).
\end{multline*}
Since  $f(0)=1$, we obtain
\[
f(c) = {(-c,-cq;q)_\infty (acd,acq/d,cdq/a,cq^2/ad;q^2)_\infty
\over (cd,cq/d,-ac,-cq/a;q)_\infty}.
\]

\noindent {\bf The sum of a very-well-poised $_6\phi_5$ series:}
\[
f(a) = {_6\phi_5} \left[ \begin{array}{c} a,qa^{\frac{1}{2}},-qa^{\frac{1}{2}},b,c,d \\
    a^{\frac{1}{2}},-a^{\frac{1}{2}},aq/b,aq/c,aq/d
    \end{array}
    ; q, {aq \over bcd} \right].
\]

By the $q$-Zeilberger algorithm, we find
\[
f(a) =
{(1-aq/cd)(1-aq/bc)(1-aq/bd)(1-aq) \over
    (1-aq/bcd)(1-aq/b)(1-aq/c)(1-aq/d)} f(aq).
\]
Since $f(0)=1$, we obtain
\[
f(a) = {(aq,aq/bc,aq/bd,aq/cd;q)_\infty \over
(aq/b,aq/c,aq/d,aq/bcd;q)_\infty}.
\]

\subsection{Two-Term Summation Formulas}

Many classical two-term nonterminating summation formulas can be
dealt with by using the same method as single summation formulas.
It turns out that for many two-term summation formulas, the two
summand share the same recurrence relation. Moreover, the boundary
values $\lim\limits_{k \to \infty} g_{0,k}$ for the two summands
cancel out.  So we still obtain homogeneous recurrence relations
which lead to infinite products.  We give three examples from the
appendix of \cite{GaRa90}, and present a detailed proof for the
first example.

\noindent {\bf 1. A nonterminating form of the $q$-Vandermonde
sum:}
\begin{multline*}
f(a,b,c) = \\
{_2\phi_1} \left[ \begin{array}{c} a,b \\
c \end{array} ; q, q \right]  + {(q/c,a,b;q)_\infty \over
(c/q,aq/c,bq/c;q)_\infty} {_2\phi_1}\left[ \begin{array}{c} aq/c, bq/c\\
q^2/c \end{array} ; q, q \right].
\end{multline*}
Since $\lim\limits_{N \to \infty} f(aq^N,bq^N,cq^N)$ does not
exist, we consider
\[
g(a,b,c) = f(a,b,c) \big/ (q/c;q)_\infty.
\]
Let
\begin{align*}
& u^{(1)}_{n,k} = {1 \over (q/cq^n;q)_\infty} {(aq^n,bq^n;q)_k
\over (cq^n,q;q)_k} q^k, \\
& u^{(2)}_{n,k} = {(aq^n,bq^n;q)_\infty \over
(cq^n/q,aq/c,bq/c;q)_\infty} {(aq/c,bq/c;q)_k \over
(q^2/cq^n,q;q)_k} q^k.
\end{align*}
We have
\[
(abq^{n+1}-c)u^{(i)}_{n+1,k} + c u^{(i)}_{n,k} =
g^{(i)}_{n,k+1}-g^{(i)}_{n,k}, \quad i=1,2,
\]
where
\begin{align*}
& g^{(1)}_{n,k} = {c (1- ab q^{2n+k}) (1-q^k) \over
  q^{k} (1-aq^n)(1-bq^n) } u^{(1)}_{n,k}, \\
& g^{(2)}_{n,k} = {c(cq^n-q^{k+1}) (1-q^k) \over
q^{k+1}(1-aq^n)(1-bq^n)} u^{(2)}_{n,k}.
\end{align*}
Noting that $g^{(1)}_{n,0}=g^{(2)}_{n,0}=0$ and
\[
\lim\limits_{k \to \infty} g^{(1)}_{n,k} = - \lim\limits_{k \to
\infty} g^{(2)}_{n,k} = - {(aq^{n+1},bq^{n+1};q)_\infty \over q^n
(1/cq^n,cq^{n+1},q;q)_\infty},
\]
we get $g(a,b,c) = (1-abq/c) g(aq,bq,cq)$. Since
\[
\lim_{N \to \infty} g(aq^N,bq^N,cq^N) =0 + {1 \over
(aq/c,bq/c;q)_\infty } = {1 \over (aq/c,bq/c;q)_\infty },
\]
we get
\begin{equation}
\label{q-Van} f(a,b,c) = (q/c;q)_\infty g(a,b,c) = {
(q/c,abq/c;q)_\infty \over (aq/c,bq/c;q)_\infty }.
\end{equation}

\noindent {\bf 2. A nonterminating form of the $q$-Saalsch\"{u}z
sum:}
\begin{multline*}
f(c) = {_3\phi_2} \left[ \begin{array}{c} a,b,c \\
e,f \end{array} ; q, q \right] + {(q/e,a,b,c,qf/e;q)_\infty \over
(e/q,aq/e,bq/e,cq/e,f;q)_\infty} \\
\cdot {_3\phi_2}\left[ \begin{array}{c} aq/e,bq/e,cq/e\\
q^2/e,qf/e \end{array} ; q, q \right],
\end{multline*}
where $f=abcq/e$.

By computation, we have
\[
f(c) = {(1-bcq/e)(1-acq/e) \over (1-cq/e)(1-abcq/e)} f(cq),
\]
and by \eqref{q-Van}
\[
\lim\limits_{N \to \infty} f(cq^N) = {(q/e,abq/e;q)_\infty \over
(aq/e,bq/e)_\infty}.
\]
Thus we get
\begin{equation}
\label{q-Saal} f(c) = {(bcq/e,acq/e,q/e,abq/e;q)_\infty \over
(cq/e,abcq/e,aq/e,bq/e;q)_\infty}.
\end{equation}

\noindent {\bf 3. Bailey's nonterminating extension of Jackson's
$_8\phi_7$ sum:}
\begin{multline*}
f(a,b) = {_8\phi_7} \left[ \begin{array}{c} a,qa^{\frac{1}{2}},-qa^{\frac{1}{2}},b,c,d,e,f \\
a^{\frac{1}{2}},-a^{\frac{1}{2}},aq/b,aq/c,aq/d,aq/e,aq/f
\end{array} ; q, q
\right] \\[10pt]
- {b \over a} {(aq,c,d,e,f,bq/a,bq/c,bq/d,bq/e,bq/f;q)_\infty
\over
(aq/b,aq/c,aq/d,aq/e,aq/f,bc/a,bd/a,be/a,bf/a,b^2q/a;q)_\infty}
\\[5pt]
\cdot {_8\phi_7}\left[ \begin{array}{c}
b^2/a,qba^{-\frac{1}{2}},-qba^{-\frac{1}{2}},b,bc/a,bd/a,be/a,bf/a \\
ba^{-\frac{1}{2}},-ba^{-\frac{1}{2}},bq/a,bq/c,bq/d,bq/e,bq/f
\end{array} ; q, q \right],
\end{multline*}
where $f=a^2q/bcde$.

By computation, we have
\[
f(a,b) = {(1-aq)(1-aq/cd)(1-aq/ce)(1-aq/de) \over
(1-aq/cde)(1-aq/c)(1-aq/d)(1-aq/e)} f(aq,bq),
\]
and by \eqref{q-Saal},
\[
\lim\limits_{N \to \infty} f(aq^N,bq^N) =
{(b/a,aq/cf,aq/df,aq/ef;q)_\infty \over
(aq/f,bc/a,bd/a,be/a;q)_\infty}.
\]
Finally, we have
\[
f(a,b) = {(aq,aq/cd,aq/ce,aq/de,b/a,aq/cf,aq/df,aq/ef;q)_\infty
\over (aq/cde,aq/c,aq/d,aq/e,aq/f,bc/a,bd/a,be/a;q)_\infty}.
\]

\subsection{Bilateral Summations}

Bilateral summations (\cite[Chapter 5]{GaRa90}) can also be dealt
with by using the $q$-Zeilberger algorithm approach. We need the
following special requirement  for the recurrence relation
\eqref{eq-ag}:
\[
\lim_{k \to -\infty} g_{n,k} = \lim_{k \to \infty} g_{n,k}=0.
\]

Here are several examples.

\noindent {\bf 1. Jacobi's triple product:}
\[
\sum_{k=-\infty}^\infty q^{k \choose 2} z^k =
(q,-z,-q/z;q)_\infty.
\]

This well-known identity is due to Jacobi \cite{Jacobi29} (see
\cite[p.~12]{Andrews71}). Cauchy \cite{Cauchy} gave a simple proof
using the $q$-binomial theorem. For other proofs, see Andrews
\cite{Andrews65}, Ewell \cite{Ewell}, Joichi and Stanton
\cite{JoSt89}.

We give a $q$-Zeilberger style proof of its semi-finite form
\cite{ChenFu05}:
\[
f(m) = \sum_{k=-\infty}^\infty {q^{k \choose 2} z^k \over
(q^{m+1};q)_k}, \quad m \ge 0.
\]
Let $u_{m,k}$ be the summand. Applying the $q$-Zeilberger
algorithm, we obtain
\[
z u_{m+1,k} - (q^{m+1}+z) (1-q^{m+1}) u_{m,k} = g_{m,k+1}-g_{m,k},
\]
where $g_{m,k} = (1-q^{m+1})q^{m+1} u_{m,k}$. Since
\[
\lim_{k \to \infty} g_{m,k} = \lim_{k \to -\infty} g_{m,k} = 0,
\]
we have
\[
f(m+1) = (1+q^{m+1}/z)(1-q^{m+1}) f(m).
\]
It follows that
\begin{multline*}
\sum_{k=-\infty}^\infty q^{k \choose 2} z^k =
\sum_{k=-\infty}^\infty \lim_{m \to \infty} u_{m,k} = \lim_{m \to \infty} f(m)\\
 = f(0) (q,-q/z;q)_\infty =
(-z,q,-q/z;q)_\infty.
\end{multline*}

\noindent
 {\bf 2. Ramanujan's $_1\psi_1$ sum:}
\[
f(b) = {_1\psi_1} \left[ \begin{array}{c}a \\b  \end{array}; q,z
\right] = {(q,b/a,az,q/az;q)_\infty \over
(b,q/a,z,b/az;q)_\infty}, \quad |z|,|b/az|<1.
\]
This formula is due to Ramanujan. Andrews
\cite{Andrews69,Andrews70}, Hahn \cite{Hahn}, Jackson
\cite{Jackson}, Ismail \cite{Ismail}, Andrews and Askey
\cite{Andrews78}, and Berndt \cite{Berndt} have found different
proofs.

The proof of Andrews and Askey \cite{Andrews78} is based on the
following recursion:
\begin{equation}
\label{1psi1-rec} f(b) = {1-b/a \over (1-b) (1-b/az)} f(bq).
\end{equation}
Instead of using series manipulations, we derive the recursion
\eqref{1psi1-rec} by using the $q$-Zeilberger algorithm. Let
$u_{n,k} = {(a;q)_k \over (bq^n;q)_k}z^k$. Then
\begin{equation}
\label{1psi1} z(bq^n-a)u_{n+1,k}+(az-bq^n)(1-bq^n)u_{n,k} =
g_{n,k+1}-g_{n,k},
\end{equation}
where
\[
g_{n,k} = (1-bq^n)bq^n \cdot u_{n,k}.
\]
Notice that \eqref{1psi1} holds for any $k \in \Z$. Furthermore,
when $|z|<1$ and $|b/az|<1$, we have
\[
\lim_{k \to \pm \infty} g_{n,k} = (1-bq^n)bq^n \cdot \lim_{k \to
\pm \infty} u_{n,k} = 0.
\]
Summing over $k \in \Z$ on both sides of \eqref{1psi1}, we
immediately get \eqref{1psi1-rec}, implying that
\[
f(b) = {(b/a;q)_\infty \over (b,b/az;q)_\infty} f(0).
\]
By the $q$-binomial theorem
\[
f(q) = \sum_{k=0}^\infty {(a;q)_k \over (q;q)_k} z^k =
{(az;q)_\infty \over (z;q)_\infty}.
\]
Therefore,
\[
f(b) = {(b/a;q)_\infty \over (b,b/az;q)_\infty} {(q,q/az;q)_\infty
\over (q/a;q)_\infty} f(q) = {(b/a,q,q/az,az;q)_\infty \over
(b,b/az,q/a,z;q)_\infty}.
\]

\noindent {\bf  3. A well-poised $_2\psi_2$ series:}
\[
f(b,c) = {_2 \psi _2} \left[ \begin{array}{c} b,c \\ aq/b,aq/c
\end{array};q,-{ aq \over bc} \right], \quad |aq/bc|<1.
\]
By computation, we have
\begin{multline*}
f(b,c) = {(1-aq/bc)(1-aq^2/bc) \over (1+aq/bc)(1+aq^2/bc)}
\\ \times {(1-aq^2/b^2)(1-aq^2/c^2) \over (1-q/b)(1-q/c)(1-aq/b)(1-aq/c)}
f(b/q,c/q).
\end{multline*}
By Jacobi's triple product identity, we obtain
\[
\lim_{N \to \infty} f(b/q^N,c/q^N) = \sum_{k=-\infty}^\infty
q^{k^2} (-a)^k = (q^2,qa,q/a;q^2)_\infty.
\]
Thus, we get
\[
f(b,c) = {(aq/bc;q)_\infty
(aq^2/b^2,aq^2/c^2,q^2,qa,q/a;q^2)_\infty \over
(-aq/bc,q/b,q/c,aq/b,aq/c;q)_\infty}.
\]

\noindent {\bf  4. Bailey's sum of a well-poised $_3\psi_3$:}
\[
f(b,c,d) = {_3 \psi _3} \left[ \begin{array}{c} b,c,d \\
q/b,q/c,q/d
\end{array};q,{q \over bcd} \right].
\]
We notice that applying the $q$-Zeilberger algorithm directly to
\[
{(b,c,d;q)_k \over (q/b,q/c,q/d;q)_k} \left({q \over
bcd}\right)^k,
\]
does not give a simple relation.  Using an idea of Paule
\cite{Paul94} of symmetrizing a bilateral summation, we replace
$k$ by $-k$ to get a summation
\[
{_3 \psi _3} \left[ \begin{array}{c} b,c,d \\
q/b,q/c,q/d
\end{array};q,{q^2 \over bcd} \right].
\]
Now we apply the $q$-Zeilberger algorithm to the average of the
above summands:
\[
{1+q^k \over 2} {(b,c,d;q)_k \over (q/b,q/c,q/d;q)_k} \left({q
\over bcd}\right)^k,
\]
and obtain that
\begin{multline*}
f(b,c,d) = {(1-q/bc)(1-q^2/bc)(1-q/bd)(1-q^2/bd) \over
(1-q/b)(1-q/c)(1-q/d)} \\
\times {(1-q/cd)(1-q^2/cd) \over (1-q/bcd)(1-q^2/bcd)(1-q^3/bcd)}
f(b/q,c/q,d/q).
\end{multline*}
By Jacobi's triple product identity, we have
\[
\lim_{N \to \infty} f(b/q^N,c/q^N,d/q^N) = \sum_{k=-\infty}^\infty
q^{3{k \choose 2}} (-q)^k = (q^3,q,q^2;q^3)_\infty = (q;q)_\infty.
\]
So we get
\[
f(b,c,d) = {(q,q/bc,q/bd,q/cd;q)_\infty \over
(q/b,q/c,q/d,q/bcd;q)_\infty}.
\]

\noindent {\bf  5. A basic bilateral analogue of Dixon's sum:}
\[
f(b,c,d) = {_4 \psi _4} \left[ \begin{array}{c} -qa,b,c,d \\
-a,a^2q/b,a^2q/c,a^2q/d
\end{array};q,{ qa^3 \over bcd} \right].
\]
By computation, we get
\begin{multline*}
f(b,c,d) = {(1-a^2q/bc)(1-a^2q^2/bc)(1-a^2q/bd)(1-a^2q^2/bd) \over
(1-a^3q/bcd)(1-a^3q^2/bcd)(1-a^3q^3/bcd)}\\
\times {(1-a^2q/cd)(1-a^2q^2/cd) \over (1-q/b)(1-q/c)(1-q/d)
} \\
\times  {(1-aq/b)(1-aq/c)(1-aq/d) \over
(1-a^2q/b)(1-a^2q/c)(1-a^2q/d)} f(b/q,c/q,d/q).
\end{multline*}
Hence,
\begin{equation}
\label{4psi4} f(b,c,d) =
{(a^2q/bc,a^2q/bd,a^2q/cd,aq/b,aq/c,aq/d;q)_\infty \over
(a^3q/bcd,q/b,q/c,q/d,a^2q/b,a^2q/c,a^2q/d;q)_\infty}  \cdot S(a),
\end{equation}
where
\[ S(a) = \lim_{N \to \infty} f(b/q^N,c/q^N,d/q^N)
= \sum_{k=-\infty}^\infty {(-qa;q)_k \over (-a;q)_k} q^{3{k
\choose 2}} (-qa^3)^k.
\]
Especially, replacing $b,c,d$ in \eqref{4psi4} by $-a,c/q^N,d/q^N$
and taking the limit $N \to \infty$, we get
\[
\lim_{N \to \infty} f(-a,c/q^N,d/q^N) = {(-q;q)_\infty \over
(-q/a,-aq;q)_\infty} \cdot S(a).
\]
By Jacobi's triple product identity, we have
\[
\lim_{N \to \infty} f(-a,c/q^N,d/q^N) = \sum_{k=-\infty}^\infty
q^{k^2} (-a^2)^k = (q^2,a^2q,q/a^2;q^2)_\infty,
\] which implies that
\[
S(a) = {(q,a^2q,q/a^2;q)_\infty \over (aq,q/a;q)_\infty}.
\]
Therefore, we obtain
\begin{multline*}
f(b,c,d) = {(a^2q/bc,a^2q/bd,a^2q/cd,aq/b,aq/c,aq/d;q)_\infty
\over
(a^3q/bcd,q/b,q/c,q/d,a^2q/b,a^2q/c,a^2q/d;q)_\infty}  \\
\times {(q,a^2q,q/a^2;q)_\infty \over (aq,q/a;q)_\infty} \qquad
\qquad
\end{multline*}

\noindent {\bf  6. Bailey's very-well-poised $_6\psi_6$ series:}
\[
f(b,c,d,e) = {_6 \psi _6} \left[ \begin{array}{c} qa^{\frac{1}{2}},-qa^{\frac{1}{2}},b,c,d,e \\
a^{\frac{1}{2}},-a^{\frac{1}{2}},aq/b,aq/c,aq/d,aq/e
\end{array};q,{ qa^2 \over bcde} \right].
\]
This identity is due to Bailey \cite{Bailey}. Other proofs have
been given by Slater and Lakin \cite{Slater56}, Andrews
\cite{Andrews74}, Chen and Liu \cite{ChenLiu}, Schlosser
\cite{Sch}, and Jouhet and Schlosser \cite{JoS}. Askey and Ismail
\cite{Askey79} gave a simple proof using $_6\phi_5$ sum and an
argument based on analytic continuation. Askey \cite{Askey84} also
showed that it can be obtained from a simple difference equation
and Ramanujan's $_1\psi_1$ sum.

Using our computational approach, we obtain
\begin{multline*}
f(b,c,d,e) = {(1-aq/bc)(1-aq^2/bc)(1-aq/bd)(1-aq^2/bd) \over
(1-aq/b)(1-aq/c)(1-aq/d)(1-aq/e)}
\\ \times {(1-aq/be)(1-aq^2/be)(1-aq/cd)(1-aq^2/cd) \over
(1-q/b)(1-q/c)(1-q/d)(1-q/e)}
 \\
\times {(1-aq/ce)(1-aq^2/ce)(1-aq/de)(1-aq^2/de) \over
(1-a^2q/bcde)(1-a^2q^2/bcde)(1-a^2q^3/bcde)(1-a^2q^4/bcde)} \\
\times f(b/q,c/q,d/q,e/q).
\end{multline*}
By Jacobi's triple product identity, we have
\begin{multline*}
\lim_{N \to \infty} f(b/q^N,c/q^N,d/q^N,e/q^N) = {1 \over 1-a}
\sum_{k=-\infty}^\infty (1-aq^{2k}) q^{4{k \choose 2}} (qa^2)^k \\
= {1 \over 1-a} \sum_{k=-\infty}^\infty  q^{k \choose 2} (-a)^k =
(q,aq,q/a;q)_\infty.
\end{multline*}
Hence we get
\[
f(b,c,d,e) =
{(aq/bc,aq/bd,aq/be,aq/cd,aq/ce,aq/de,q,aq,q/a;q)_\infty \over
(aq/b,aq/c,aq/d,aq/e,q/b,q/c,q/d,q/e,qa^2/bcde;q)_\infty}.
\]

\section{Transformation Formulas} \label{S3}

In this section, we show that many classical transformation
formulas of nonterminating basic hypergeometric series can be
proved by using the $q$-Zeilberger algorithm. The basic idea is to
find the same recurrence relation and limit value of two
summations $f(a, \ldots, c)$ and $g(a, \ldots, c)$. Suppose we
have obtained a recurrence relation of second order or higher
order of the form \eqref{rec} for both $f(a, \ldots, c)$ and $g(a,
\ldots, c)$, then
 the following theorem ensures that
$f(a,\ldots,c)$ and $g(a, \ldots, c)$ must be equal as long as
$\lim\limits_{N \to \infty} f(aq^N,\ldots,cq^N)$ coincides with
the limit $\lim\limits_{N \to \infty} g(aq^N,\ldots,cq^N)$.

\begin{theo} \label{th-rec}
Let $f(z)$ be a continuous function defined on the disc $|z| \le
r$ and $d \ge 2$ be an integer. Suppose that we have a recurrence
relation
\begin{equation}
\label{eq-rec} f(z)= a_1(z)f(zq) + a_2(z) f(zq^2) + \cdots +
a_d(z)f(zq^d).
\end{equation}
For $i=1,\ldots,d$, we denote $a_i(0)$ by $w_i$. Suppose that
there exists a real number $M>0$ such that
\[
|a_i(z) - w_i| \le M |z|, \quad 1 \le i \le d,
\]
and
\begin{align*}
& |w_d| + |w_{d-1}+w_d| + \cdots + |w_2+\cdots+w_d| < 1, \\[5pt]
& w_1+w_2+\cdots+w_d = 1.
\end{align*}
Then $f(z)$ is uniquely determined by $f(0)$ and the functions
$a_i(z)$.
\end{theo}
\pf By the recurrence relation \eqref{eq-rec}, we have
\[
f(z) = \sum_{i=1}^d A_n^{(i)} f(zq^{n+i}),
\]
where $A^{(i)}_0=a_i(z)$ and
\begin{equation}
\label{eq-Arec} \begin{cases} A^{(i)}_{n+1} = a_i(zq^{n+1})
A^{(1)}_n + A^{(i+1)}_n, & 1 \le i < d, \\[8pt]
A^{(d)}_{n+1} = a_d(zq^{n+1}) A^{(1)}_n.
\end{cases}
\end{equation}

Let
\[
\lambda(x) = x^{d-1} - \sum_{i=2}^d  \left| \sum_{j=i}^d w_j
\right| x^{d-i}.
\]
By the assumption, $\lambda(1)>0$. Hence we may choose a real
number $p$ such that $|q| < p <1$ and $\lambda(p)>0$, namely,
\[
 \sum_{i=2}^d p^{d-i} \left| \sum_{j=i}^d w_j \right| < p^{d-1}.
\]
 Let
\begin{align*}
& A=\max\left\{ |A^{(1)}_0|,\ldots, |A^{(1)}_d| \right\}, \\
& A' = \max\left\{ A \cdot dMr/\lambda(p),\,
|A^{(1)}_1-A^{(1)}_0|/p, \, \ldots,
  |A^{(1)}_d-A^{(1)}_{d-1}|/p^d  \right\}, \\
& B=dMr/p^{d-2}+A'p/A.
\end{align*}
We will use induction on $n$ to show that
\begin{eqnarray}
& |A^{(1)}_n| \le A \cdot (-B;p)_n, \label{Bound-A}\\[8pt]
& |A^{(1)}_n-A^{(1)}_{n-1}| \le A' \cdot p^n \cdot (-B;p)_n.
\label{Bound-DiffA}
\end{eqnarray}
By definition, the inequalities \eqref{Bound-A} and
\eqref{Bound-DiffA} hold for $n=1,\ldots,d$. Suppose $n \ge d$ and
the inequalities hold for $1,2,\ldots,n$. From \eqref{eq-Arec} it
follows that
\begin{eqnarray*}
|A^{(1)}_{n+1}| &=& \big| \sum_{i=1}^d a_i(zq^{n+2-i}) A^{(1)}_{n+1-i}  \big| \\
& \le & \sum_{i=1}^d \left( \left| a_i(zq^{n+2-i})-w_i \right|
\cdot |A^{(1)}_{n+1-i}| \right) +
|A^{(1)}_n| \\
&& \quad +  \left| \sum_{i=2}^d \left(
(A^{(1)}_{n+1-i}-A^{(1)}_{n+2-i}) \sum_{j=i}^d w_j  \right)
\right|
\end{eqnarray*}
By the inductive hypotheses, it follows that
\begin{eqnarray*}
 |A^{(1)}_{n+1}| & \le & \sum_{i=1}^d Mr|q^{n+2-i}| \cdot A \cdot (-B;p)_{n+1-i} + A \cdot
 (-B;p)_n \\
&& \quad + A' \sum_{i=2}^d p^{n+2-i} (-B;p)_{n+2-i} \left|
\sum_{j=i}^d w_j \right| \\
& \le & A (1+dMr/p^{d-2}\cdot p^n) (-B;p)_n  + A' \cdot (-B;p)_n
\sum_{i=2}^d p^{n+2-i} \left|
\sum_{j=i}^d w_j \right|\\
& < & A \big( 1+(dMr/p^{d-2}+A'p/A)\cdot p^n \big) (-B;p)_n \\
&=& A \cdot (-B;p)_{n+1}.
\end{eqnarray*}
Similarly, by the inductive assumptions we have
\begin{eqnarray*}
\lefteqn{|A^{(1)}_{n+1}-A^{(1)}_n|}  \\
& \le & \sum_{i=1}^d \left| a_i(zq^{n+2-i})-w_i \right|
  |A^{(1)}_{n+1-i}|  +  \left| \sum_{i=2}^d \left(
  (A^{(1)}_{n+1-i}-A^{(1)}_{n+2-i})
  \sum_{j=i}^d w_j  \right) \right| \\
& \le & A \cdot dMr/p^{d-1}\cdot p^{n+1} (-B;p)_n  + A' \cdot
(-B;p)_n \sum_{i=2}^d p^{n+2-i} \left|
\sum_{j=i}^d w_j \right|\\
& = & A' \big( AdMr/(A'p^{d-1}) + 1 - \lambda(p)/p^{d-1})\cdot
p^{n+1} (-B;p)_n
\\
& \le & A' p^{n+1}  (-B;p)_{n+1}.
\end{eqnarray*}
Therefore, the inequalities  \eqref{Bound-A} and
\eqref{Bound-DiffA} hold for $n+1$. Using \eqref{Bound-DiffA} we
reach  the following inequality
\[
|A^{(1)}_n-A^{(1)}_{n-1}| \le  A' p^n  (-B;p)_\infty .
\]
So the limit $\lim\limits_{n \to \infty} A^{(1)}_n$ exists. By
\eqref{eq-Arec}, for any $1 \le i \le d$, the $\lim\limits_{n \to
\infty} A^{(i)}_n$ exists. Thus,  we get
\[
f(z) = f(0) \sum_{i=1}^d \lim\limits_{n \to \infty}  A^{(i)}_n ,
\]
which completes the proof. \qed

\noindent {\bf Remarks.}
\begin{itemize}
\item The condition that $f(z)$ is continuous in $|z| \le r$ can
be replaced by the assumption that $\lim_{N \to \infty} f(zq^N)$
exists.

\item The above theorem can be easily generalized to
multi-variable.
\end{itemize}

We now give some examples. The first five examples are adopted
from the appendix of \cite{GaRa90}.

\noindent {\bf 1.  Heine's transformations of $_2\phi_1$ series:}
\begin{eqnarray}
{_2 \phi _1} \left[ \begin{array}{c} a,b \\ c \end{array}
;q,z\right] &=& {(b,az;q)_\infty \over (c,z;q)_\infty}
{_2 \phi _1} \left[ \begin{array}{c} c/b,z \label{Heine1}\\
az
\end{array} ;q,b\right] \\
&=& {(c/b,bz;q)_\infty \over (c,z;q)_\infty}
{_2 \phi _1} \left[ \begin{array}{c} abz/c,b \\
bz
\end{array} ;q,c/b\right] \label{Heine2}\\
&=& {(abz/c;q)_\infty \over (z;q)_\infty}
{_2 \phi _1} \left[ \begin{array}{c} c/a,c/b \\
c
\end{array} ;q,abz/c\right]. \label{Heine3}
\end{eqnarray}

Let
\[
f(z) = {_2 \phi _1} \left[ \begin{array}{c} a,b \\ c \end{array}
;q,z\right].
\]
We have
\begin{equation} \label{eqHeine} f(z) = {-c-q+(qa+qb)z \over
q(z-1)} f(zq) + {c-qabz \over q(z-1)} f(zq^2).
\end{equation}
By Theorem~\ref{th-rec}, for $|c/q|<1$, $f(z)$ is uniquely
determined by $f(0)$ and the recurrence relation \eqref{eqHeine}.
Let
\[
g(z) = {(b,az;q)_\infty \over (c,z;q)_\infty}
{_2 \phi _1} \left[ \begin{array}{c} c/b,z \\
az
\end{array} ;q,b\right].
\]
Then $g(z)$ satisfies the same recursion as \eqref{eqHeine}. By
the $q$-binomial theorem, we have
\[
g(0) =  {(b;q)_\infty \over (c;q)_\infty}
{_1 \phi _0} \left[ \begin{array}{c} c/b \\
\text{---}
\end{array} ;q,b\right] = 1 = f(0).
\]
Therefore, \eqref{Heine1} holds for $|c/q|<1$. By analytic
continuation, \eqref{Heine1} holds for all $a,b,c,z \in \C$
provided that both sides are convergent. Similar arguments can
justify \eqref{Heine2} and \eqref{Heine3}.

\noindent {\bf 2. Jackson's transformations of $_2\phi_1$,
$_2\phi_2$ and
 $_3\phi_2$ series:}
\begin{eqnarray}
{_2 \phi _1} \left[ \begin{array}{c} a,b \\ c \end{array}
;q,z\right] &=& {(az;q)_\infty \over (z;q)_\infty}
{_2 \phi _2} \left[ \begin{array}{c} a, c/b \\
c,az
\end{array} ;q,bz\right] \label{2p2}\\
&=& {(abz/c;q)_\infty \over (bz/c;q)_\infty}
{_3 \phi _2} \left[ \begin{array}{c} a,c/b,0 \\
c,cq/bz
\end{array} ;q,q\right], \label{3p2}
\end{eqnarray}
where \eqref{3p2} holds provided that the series terminates.

Let $f(z)$ be the left hand side \eqref{2p2}. Thus we have the
recurrence relation \eqref{eqHeine} and $\lim\limits_{N \to
\infty}f(zq^N)=1$. By using the $q$-Zeilberger algorithm, one can
verify that the right hand sides of \eqref{2p2} and \eqref{3p2}
also satisfy the same recurrence relation. Moreover, for the
summation \eqref{3p2} the terminating condition is required to
ensure $\lim\limits_{k \to \infty} g_{n,k}=0$ in \eqref{eq-ag}. By
considering the limit values, we get the transformation formulas
\eqref{2p2} and \eqref{3p2}.

A similar discussion implies the following transformation formula
for terminating $_2\phi_1$ series:
\[
{_2 \phi _1} \left[ \begin{array}{c} a,b \\ c \end{array}
;q,z\right] = {(c/b,c/a;q)_\infty \over (c/ab,c;q)_\infty}
{_3 \phi _2} \left[ \begin{array}{c} a,b,abz/c \\
abq/c,0
\end{array} ;q,q\right],
\]
provided that the right hand side summation terminates.

\noindent {\bf 3.  Transformations of $_3\phi_2$ series:}
\begin{eqnarray}
\lefteqn{{_3 \phi _2} \left[ \begin{array}{c} a,b,c \\ d,e
\end{array} ;q,{de \over abc}\right]} \label{3p2deabc} \\
&=& {(e/a,de/bc;q)_\infty \over (e,de/abc;q)_\infty}
{_3 \phi _2} \left[ \begin{array}{c} a,d/b,d/c \\
d,de/bc
\end{array} ;q,{e \over a}\right] \label{III.9}\\
&=& {(b,de/ab,de/bc;q)_\infty \over (d,e,de/abc;q)_\infty}
{_3 \phi _2} \left[ \begin{array}{c} d/b,e/b,de/abc \\
de/ab,de/bc
\end{array} ;q,b\right]. \label{III.10}
\end{eqnarray}

We take $d$ as the parameter.  Let $f(d)$ be the series in
\eqref{3p2deabc}. We have $f(0)=1$ and
\begin{multline*}
f(d) =  -{(1+q)ed\,^2+(-abc-eb-ea-ec)d+abc+abce/q \over
(-ed+abc)(-1+d)} f(dq) \\
+ {e(-c+dq)(-dq+b)(-dq+a) \over q(-ed+abc)(-1+dq)(-1+d)} f(dq^2).
\end{multline*}
On the other hand, one can verify that both the series in
\eqref{III.9} and \eqref{III.10} have the same limit value and
satisfy the same
 recurrence relation as \eqref{3p2deabc}.

\noindent {\bf 4. Sears-Carlitz transformation:}
\begin{multline*}
{_3\phi_2} \left[ \begin{array}{c} a,\ b,\ c \\
aq/b,\ aq/c \end{array}; q, {aqz \over bc} \right] \\
= {(az;q)_\infty \over (z;q)_\infty} {_5\phi_4} \left[
\begin{array}{c} a^{\frac{1}{2}},\ -a^{\frac{1}{2}},
\ (aq)^{\frac{1}{2}},\ -(aq)^{\frac{1}{2}},
\ aq/bc \\
aq/b,\ aq/c,\ az,\ q/z \end{array}; q, q \right],
\end{multline*}
provided that the right hand side terminates.

Let us take $z$ as the parameter and denote the series by $f(z)$.
One can verify that both sides have the same limit value
$\lim\limits_{N \to \infty} f(zq^N) = 1$ and satisfy the following
recurrence relation:
\[
f(z)=r_1(z)f(zq)+r_2(z)f(zq^2)+r_3(z)f(zq^3),
\]
where
\begin{align*}
&r_1(z)= {ab+ac+bc \over bc}+O(z), \quad r_2(z) = {-a(b+a+c) \over
bc}+O(z), \\
&r_3(z)={a^2 \over bc}+O(z).
\end{align*}
Note that to comply with the conditions of Theorem~\ref{th-rec},
we only need the values  $r_1(0)$, $r_2(0)$ and $r_3(0)$. So we do
not give the explicit formulas for $r_1(z)$, $r_2(z)$ and
$r_3(z)$.

 \noindent {\bf 5. Transformations of very-well-poised
$_8\phi_7$ series:}

\begin{eqnarray}
\lefteqn{{_8\phi_7} \left[ \begin{array}{c} a, qa^{\frac{1}{2}},-qa^{\frac{1}{2}},b,c,d,e,f \\
a^{\frac{1}{2}},-a^{\frac{1}{2}},aq/b,aq/c,aq/d,aq/e,aq/f
\end{array}; q, {a^2q^2 \over bcdef}
\right]} \label{8p7aa}   \\
&=& {(aq,aq/ef,\lambda q/e,\lambda q/f;q)_\infty \over
(aq/e,aq/f,\lambda q, \lambda q/ef;q)_\infty} \nonumber \\
&& \times {_8\phi_7} \left[ \begin{array}{c} \lambda, q
\lambda^{\frac{1}{2}}, -q \lambda^{\frac{1}{2}},
\lambda b/a, \lambda c/a, \lambda d/a, e,f\\
\lambda^{\frac{1}{2}}, -\lambda^{\frac{1}{2}}, aq/b, aq/c, aq/d,
\lambda q/e, \lambda q/f
\end{array}; q, {aq \over ef} \right] \label{8p7-1}\\
&=& {(aq,b,bc\mu/a,bd\mu/a,be\mu/a,bf\mu/a;q)_\infty \over
(aq/c,aq/d,aq/e,aq/f,\mu q,b\mu/a;q)_\infty} \nonumber \\
&& \times  {_8\phi_7} \left[ \begin{array}{c} \mu,
q\mu^{\frac{1}{2}}, -q\mu^{\frac{1}{2}},
aq/bc,aq/bd,aq/be,aq/bf,b\mu/a\\
\mu^{\frac{1}{2}},-\mu^{\frac{1}{2}},bc\mu/a,bd\mu/a,be\mu/a,bf\mu/a,aq/b
\end{array}; q, b \right], \nonumber \\ \label{8p7-2}
\end{eqnarray}
where $\lambda=qa^2/bcd$ and $\mu=q^2a^3/b^2cdef$.

 We choose  $a,b,f$ as parameters for the series in
 \eqref{8p7aa} and \eqref{8p7-1} and denote the series by $H(a,b,f)$.
 It follows from \eqref{III.9}
 that  both series have the same
 limit value  $\lim\limits_{N
\to \infty}H(aq^N,bq^N,fq^N)$. By computation, one sees that they
satisfy the following recurrence relation
\[
H(a,b,f)=r_1(a,b,f)H(aq,bq,fq)+r_2(a,b,f)H(aq^2,bq^2,fq^2),
\]
where
\[
r_1(a,b,f)= 1+O(a), \quad r_2(a,b,f) = O(a).
\]
Thus, we have verified the first transformation formula. To prove
the second transformation formula, we choose $a,c,f$ as the
parameters and denote the series by $H(a,c,f)$. By computation,
the series in \eqref{8p7aa} and \eqref{8p7-2} satisfy the
following recurrence relation
\[
H(a,c,f)=r_1(a,c,f)H(aq,cq,fq)+r_2(a,c,f)H(aq^2,cq^2,fq^2).
\]
 Using the transformation formula \eqref{III.10}, one sees that
 both sides have the same limit value $\lim\limits_{N
\to \infty}H(aq^N,cq^N,fq^N)$. Thus we have obtained the second
transformation formula.

\noindent {\bf 6. A limiting case of Watson's formula.}

Watson \cite{Watson} used the following formula to prove the
Rogers-Ramanujan identities \cite{Hardy} (see also \cite[Section
2.7]{GaRa90}):
\begin{equation}
\label{RR} \sum_{k=0}^\infty {(aq;q)_{k-1} (1-aq^{2k}) \over
(q;q)_k} (-1)^ka^{2k}q^{k(5k-1)/2} = (aq;q)_\infty
\sum_{k=0}^\infty {a^k q^{k^2} \over (q;q)_k}.
\end{equation}

We choose $a$ as the parameter. Then we can verify that both sides
of \eqref{RR} have the same limit value $f(0)=1$ and satisfy the
same recurrence relation
\[ f(a)= (1-aq)f(aq) + aq(1-aq)(1-aq^2)
f(aq^2).
\]

Setting $a=1$ and $a=q$ in \eqref{RR}, we obtain the
Rogers-Ramanujan identities by Jacobi's triple product identity:
\begin{eqnarray*}
(q;q)_\infty \sum_{k=0}^\infty {q^{k^2} \over (q;q)_k} &=&
\sum_{k=-\infty}^\infty (-q^2)^k q^{5{k \choose 2}} =
(q^2,q^3,q^5;q^5)_\infty,
\end{eqnarray*}
and
\begin{eqnarray*}
(q;q)_\infty \sum_{k=0}^\infty {q^{k^2+k} \over (q;q)_k} &=&
\sum_{k=-\infty}^\infty (-q^4)^k q^{5{k \choose 2}} =
(q,q^4,q^5;q^5)_\infty.
\end{eqnarray*}
Finite forms of the above identities have been proved by Paule
\cite{Paul94} by using the $q$-Zeilberger algorithm.

\noindent {\bf 7. A generalization of Lebesgue's identity.}

The following transformation formula is due to Carlitz
\cite{Car62} (see \cite{Andrews66}):
\begin{equation}\label{g-Leb}
\sum_{k=0}^\infty {(x;q)_k q^{k \choose 2} (-a)^k \over
(q,bx;q)_k} = {(a,x;q)_\infty \over (bx;q)_\infty}
\sum_{k=0}^\infty {(b;q)_k x^k \over (q,a;q)_k}.
\end{equation}

We choose $x$ as the parameter. Both sides of \eqref{g-Leb} have
the same limit value $f(0)=(a;q)_\infty$ and satisfy the same
recurrence relation:
\[ f(x) = \left({q+a \over q} +
O(x)\right)f(xq) + \left({-a \over q} + O(x)\right) f(xq^2).
\]

\noindent {\bf 8. Three-Term transformation formulas.}

 Our approach also applies to certain three-term transformation formulas.
 It is sometimes the case that the left hand side of the
identity satisfies a homogenous recursion, and the two terms on
the right hand side satisfy non-homogenous recursions respectively
but their sum leads to a homogenous recurrence relation.

The first  example is
\begin{eqnarray}
\lefteqn{{_3\phi_2} \left[ \begin{array}{c} a,b,c \\
d,e \end{array}; q, {de \over abc} \right]} \nonumber \\
&=& {(e/b,e/c;q)_\infty \over (e,e/bc;q)_\infty}
\nonumber {_3\phi_2} \left[ \begin{array}{c} d/a,b,c \\
d,bcq/e \end{array}; q, q \right] \\
&&+ {(d/a,b,c,de/bc;q)_\infty \over (d,e,bc/e,de/abc;q)_\infty}
 {_3\phi_2} \left[
\begin{array}{c} e/b,e/c,de/abc \\
de/bc,eq/bc
\end{array}; q, q \right]. \label{III.34}
\end{eqnarray}

Let us choose $e$ as the parameter. Then both sides of
\eqref{III.34} have the same limit value $\lim\limits_{N \to
\infty} f(eq^N) = 1$ and  satisfy the same recurrence relation:
\[
f(e) = r_1(e) f(eq) + r_2(e) f(eq^2),
\]
where
\[
r_1(e) = {q+d \over q} + O(e), \quad r_2(e) = -{d \over q} + O(e).
\]

The second example  is
\begin{eqnarray}
\lefteqn{{_8\phi_7} \left[ \begin{array}{c} a,qa^{\frac{1}{2}},-qa^{\frac{1}{2}},b,c,d,e,f \\
a^{\frac{1}{2}},-a^{\frac{1}{2}},aq/b,aq/c,aq/d,aq/e,aq/f
\end{array}; q, {a^2q^2 \over bcdef}
\right]} \nonumber \\
&=& {(aq,aq/de,aq/df,aq/ef;q)_\infty \over
(aq/d,aq/e,aq/f,aq/def;q)_\infty}
\nonumber {_4\phi_3} \left[ \begin{array}{c} aq/bc,d,e,f \\
aq/b,aq/c,def/a \end{array}; q, q \right] \\
&&+ {(aq,aq/bc,d,e,f,a^2q^2/bdef,a^2q^2/cdef;q)_\infty \over
(aq/b,aq/c,aq/d,aq/e,aq/f,a^2q^2/bcdef,def/aq;q)_\infty} \nonumber \\
&& \qquad \cdot {_4\phi_3} \left[
\begin{array}{c} aq/de,aq/df,aq/ef,a^2q^2/bcdef \\
a^2q^2/bdef,a^2q^2/cdef,aq^2/def
\end{array}; q, q \right]. \label{III.36}
\end{eqnarray}

We take $a,b,f$ as parameters and denote the series by $H(a,b,f)$.
By the transformation formula \eqref{III.34}, we see that both
sides of \eqref{III.36} have the same limit value $\lim\limits_{N
\to \infty} H(aq^N,bq^N,fq^N)$. Moreover, they satisfy the same
recurrence relation:
\[ H(a,b,f) = r_1(a,b,f) H(aq,bq,fq) +
r_2(a,b,f) H(aq^2,bq^2,fq^2),
\]
where
\[
r_1(a,b,f) = 1 + O(a), \quad r_2(a,b,f) =  O(a).
\]

\noindent {\bf 9. The Rogers-Fine identity.}

To conclude this paper, we consider a transformation formula that
can be justified by using nonhomogeneous recurrence relations.
This is the Rogers-Fine identity \cite{Fine}:
\begin{equation}\label{Rogers-Fine}
\sum_{k=0}^\infty {(a;q)_k \over (b;q)_k} z^k = \sum_{k=0}^\infty
{(a,azq/b;q)_k (1-azq^{2k}) q^{k^2-k} (bz)^k \over
(b,z;q)_k(1-zq^k)}.
\end{equation}
We choose  $z$ as the parameter.  By computation, both sides of
\eqref{Rogers-Fine} satisfy the following recurrence relation:
\[
(-b+azq) f(zq) + (-zq+q) f(z) = q-b.
\]
The non-homogenous term $q-b$ appears because $g_{0,0}=-q+b$ and
$\lim\limits_{k \to \infty} g_{0,k}=0$ when one implements the
$q$-Zeilberger algorithm. Let $d(z)$ be the difference of the two
sides of \eqref{Rogers-Fine}. Then we have
\[
d(z) = {b \over q} \cdot {1-azq/b \over 1-z} d(zq).
\]
Since $d(0)=0$, the equation \eqref{Rogers-Fine} holds for $|z|<1$
and $|b|<|q|$. By analytic continuation, it holds for $|z|<1$.

\vskip 15pt \noindent {\small {\bf Acknowledgments.} This work was
supported by the ``973'' Project on Mathematical Mechanization,
the National Science Foundation, the Ministry of Education, and
the Ministry of Science and Technology of China.}

\end{document}